# New variances for various kappa coefficients based on the unbiased estimator of the expected index of agreements

*Martín Andrés, A.*[1]* *and Álvarez Hernández, M.*[2,3]

[1] Department of Biostatistics, Faculty of Medicine, University of Granada, Spain.
[2] CITMAga, 15782 Santiago de Compostela, Spain.
[3] Defense University Center, Spanish Naval Academy, Marín, Pontevedra, Spain.

## ABSTRACT

Recently Martín Andrés and Álvarez Hernández (2024) have proposed new estimators of various kappa coefficients. These estimators are based on the unbiased estimator of the expected index of agreement of each population coefficient. In their article, these authors propose variance formulas based on the univariate delta method. Here new formulas are proposed that are based on the multivariate delta method.

**Keywords:** Agreement; Cohen's *kappa*; concordance and intraclass correlation coefficients; Conger's *kappa*; Fleiss' *kappa*; Gwet's *AC1/2*; Hubert's *kappa*; Krippendorf's *alpha*; pairwise multi-rater *kappa*; Scott's *pi*.

## 1.- Obtaining the new variances.

For everything that follows, the same notation and numbering of formulas as defined in the cited article [1] is used, since what is current should be understood as a complement to Appendix 2 of said article. Let us assume that the variance V($\hat{\kappa}_X$) of a classic estimator $\hat{\kappa}_X$ has been obtained (or can be obtained) through the multivariant version of the delta method. Then another expression $V_A(\hat{\kappa}_{XU})$ can be deduced from the variance of the new estimator $\hat{\kappa}_{XU}$, which is different to the variance V($\hat{\kappa}_{XU}$) in Appendix 2.

* Correspondence to: Bioestadística, Facultad de Medicina, C8-01, Universidad de Granada, 18071 Granada, Spain. Email: amartina@ugr.es. Phone: 34-58-244080.

Let this be the case of $R$=2. If $\boldsymbol{x}=\{x_{ij}\}$ is a set of observations from a multinomial distribution of parameters $\boldsymbol{p}=\{p_{ij}\}$ and $n$, let $\hat{p}_{ij}=x_{ij}/n$ and $\hat{\boldsymbol{p}}=\{\hat{p}_{ij}\}$ be the unbiased estimations of these parameters, and let $f(\boldsymbol{p})$ be any function, estimated by $\hat{f}=f(\hat{\boldsymbol{p}})$, with derivatives $f_{ij}=\partial f/\partial p_{ij}$. The multivariant delta method indicates that $V(\hat{f})=[\Sigma\Sigma f_{ij}^2 p_{ij}-(\Sigma\Sigma f_{ij}p_{ij})^2]/n$ [2]. When $\hat{f}=\hat{\kappa}_C$, then $V(\hat{\kappa}_C)$ is given by the classic expression of Fleiss *et al.* [3], which can be deduced through the multivariant version of the delta method [4]. If $\hat{g}=\hat{\kappa}_{CU}$ then, through expression (4), $g_{ij}=(\partial g/\partial\kappa_C)\times(\partial f/\partial p_{ij})=(\partial g/\partial\kappa_C)\times f_{ij}=[n(n-1)/(n-1+\kappa_C)^2]\times f_{ij}$. Therefore $V(\hat{\kappa}_{CU})=[n(n-1)/(n-1+\kappa_C)^2]^2\times[\Sigma\Sigma f_{ij}^2 p_{ij}-(\Sigma\Sigma f_{ij}p_{ij})^2]/n$ and thus the new variance -alternative to expression (5)- is

$$V_A\left(\hat{\kappa}_{CU}\right)=\frac{\left\{n\left(n-1\right)\right\}^2}{\left(n-1+\kappa_C\right)^4}V\left(\hat{\kappa}_C\right). \qquad \textbf{(5bis)}$$

Something similar occurs in the following cases. If $\hat{g}=\hat{\kappa}_{SU}$ then, through expression (12), $\partial\kappa_{SU}/\partial\kappa_S=4n(n-1)/(2n-1+\kappa_S)^2$ and therefore the alternative to expression (13) is

$$V_A\left(\hat{\kappa}_{SU}\right)=\frac{\left\{4n\left(n-1\right)\right\}^2}{\left(2n-1+\kappa_S\right)^4}V\left(\hat{\kappa}_S\right). \qquad \textbf{(13bis)}$$

The case of $\hat{\kappa}_{KU}$ is deduced through the traditional delta method. Through expression (15), $\partial\kappa_{KU}/\partial\kappa_{SU}=(2n-1)/2n$, $V(\hat{\kappa}_{KU})=[(2n-1)/2n]^2 V(\hat{\kappa}_{SU})$ and, through expression (13bis), we obtain the following alternative to the variance of the end of Section 2.4:

$$V_A\left(\hat{\kappa}_{KU}\right)=\frac{\left\{2\left(2n-1\right)\left(n-1\right)\right\}^2}{\left(2n-1+\kappa_S\right)^4}V\left(\hat{\kappa}_S\right).$$

Let this now be the case of $R$>2, in which the expressions at the beginning of Paragraph 2 are valid, with the necessary changes. As indicated in Section 3.2, the variance for $\hat{\kappa}_{HU}$ is like that for $\hat{\kappa}_{CU}$, changing the letter $C$ to the letter $H$; therefore

$$V_A\left(\hat{\kappa}_{HU}\right)=\frac{\left\{n\left(n-1\right)\right\}^2}{\left(n-1+\kappa_H\right)^4}V\left(\hat{\kappa}_H\right).$$

If $\hat{g}=\hat{\kappa}_{FU}$ then, through expression (33), $\partial\kappa_{FU}/\partial\kappa_F=R^2n(n-1)/\{R(n-1)+1+(R-1)\kappa_F\}^2$ and therefore an alternative to expression (34) is

$$V_A\left(\hat{\kappa}_{FU}\right)=\frac{\left\{R^2n\left(n-1\right)\right\}^2}{\left\{nR-\left(R-1\right)\left(1-\kappa_F\right)\right\}^4}V\left(\hat{\kappa}_F\right). \qquad \textbf{(34bis)}$$

If $\hat{g}=\hat{\kappa}_{F2U}$, through the end of Section 3.3 it is obtained that $\partial\kappa_{F2U}/\partial\kappa_{F2}=4n(n-1)/(2n-1+\kappa_{F2})^2$; therefore, the alternative variance to that of the end of the aforementioned section is

$$V_A\left(\hat{\kappa}_{F2U}\right)=\frac{\left\{4n\left(n-1\right)\right\}^2}{\left(2n-1+\kappa_{F2}\right)^4}V\left(\hat{\kappa}_{F2}\right).$$

Finally, through what is indicated at the end of Section 3.4, variance $V_A(\hat{\kappa}_{KU})$ is that of the case $R$=2 changing the letter $S$ to the letter $F$; therefore, the alternative to expression (38) is

$$V_A\left(\hat{\kappa}_{KU}\right)=\frac{\left\{2\left(2n-1\right)\left(n-1\right)\right\}^2}{\left(2n-1+\kappa_F\right)^4}V\left(\hat{\kappa}_F\right), \qquad \textbf{(38bis)}$$

as the variance of $V_A(\hat{\kappa}_{K2U})$ is the same as above, but adding the number 2 to the letters $K$ and $F$

$$V_A\left(\hat{\kappa}_{K2U}\right)=\frac{\left\{2\left(2n-1\right)\left(n-1\right)\right\}^2}{\left(2n-1+\kappa_{F2}\right)^4}V\left(\hat{\kappa}_{F2}\right)$$

An interesting observation is that the values of $\kappa_X$ in which $V_A(\hat{\kappa}_{XU})=V(\hat{\kappa}_X)$ are also the values of $\kappa_X$ where the difference $\hat{\kappa}_{XU}-\hat{\kappa}_X$ becomes extreme.

## 2.- Simulation to assess the two variances of $\hat{\kappa}_{CU}$ (unweighted).

For this objective, it is necessary to proceed in the same way as in Section 5 of the article, with the appropriate changes. The steps are:

(1) From the multinomial distribution $\{p_{ij}; n\}$ that is chosen, all of its populational parameters are known and, in particular, its real *kappa* ($\kappa_C$) value. From these we can determine the exact asymptotic variances $V(\hat{\kappa}_{CU})$ and $V_A(\hat{\kappa}_{CU})$ from expressions (5) and (5bis) respectively, where [3]

$$V\left(\hat{\kappa}_C\right)=\frac{A+B-C}{n\left(1-I_e\right)^2} \text{ with } \begin{cases} A=\sum_i p_{ii}\left[1-\left(p_{i\bullet}+p_{\bullet i}\right)\left(1-\kappa_C\right)\right]^2, \\ B=\left(1-\kappa_C\right)^2\sum_i\sum_{j\neq i} p_{ij}\left(p_{\bullet i}+p_{j\bullet}\right)^2, \\ C=\left[\kappa_C-I_e\left(1-\kappa_C\right)\right]^2=\left[1-\left(1-\kappa_C\right)\left(1+I_e\right)\right]^2. \end{cases}$$

(2) Extract $N$=10.000 samples of that multinomial distribution, calculating for each one of them ($h$ = 1, 2,…, $N$) the values of:

→ The estimation $\hat{\kappa}_{CUh}$ of expression (4) of the article.

→ The estimated variances $\hat{V}_h\left(\hat{\kappa}_{CU}\right)$ and $\hat{V}_{Ah}\left(\hat{\kappa}_{CU}\right)$ of the exact asymptotic variances $V(\hat{\kappa}_{CU})$ and $V_A(\hat{\kappa}_{CU})$, which are obtained by substituting the populational values $\{p_{i\bullet}, p_{\bullet j}, p_{ij}, I_e, \kappa_C\}$ with the sample values $\left\{\hat{p}_{i\bullet}, \hat{p}_{\bullet j}, \hat{p}_{ij}, \hat{I}_{eU}, \hat{\kappa}_{CU}\right\}$ of sample $h$.

(3) With this data we determine:

→ Variance $\hat{V}_E\left(\hat{\kappa}_{CU}\right)$ (with denominator $N–1$) of the $N$ values $\hat{\kappa}_{CUh}$, which is an estimation of the exact variance $V_E\left(\hat{\kappa}_{CU}\right)$. As the latter is unknown and $N$ is large enough, we can understand that $\hat{V}_E\left(\hat{\kappa}_{CU}\right)\approx V_E\left(\hat{\kappa}_{CU}\right)$.

→ The average value $\bar{\hat{V}}\left(\hat{\kappa}_{CU}\right)$ and $\bar{\hat{V}}_A\left(\hat{\kappa}_{CU}\right)$ of the $N$ estimators $\hat{V}_h\left(\hat{\kappa}_{CU}\right)$ and $\hat{V}_{Ah}\left(\hat{\kappa}_{CU}\right)$.

(4) The process is repeated for scenarios $K$ = 2, 3, and 5; $n$ = 10, 20, 50, and 100; $\kappa_C \approx 0.4$ and 0.8. The results obtained are those for Table 6 (attached).

It can be observed that $V(\hat{\kappa}_{CU})$ and $V_A(\hat{\kappa}_{CU})$ usually give values lower than $\hat{V}_E\left(\hat{\kappa}_{CU}\right)$,

with a relative average of −7.6% and −7.9% respectively, and $V(\hat{\kappa}_{CU})$ is somewhat better. A similar thing happens with $\hat{\bar{\mathrm{V}}}(\hat{\kappa}_{CU})$ and $\hat{\bar{\mathrm{V}}}_A(\hat{\kappa}_{CU})$ in relation to $\hat{\mathrm{V}}_E(\hat{\kappa}_{CU})$; now the relative averages are −7.2% and −7.5% respectively, so that $\hat{\bar{\mathrm{V}}}(\hat{\kappa}_{CU})$ is somewhat better. Therefore, the estimator selected is $\hat{\bar{\mathrm{V}}}(\hat{\kappa}_{CU})$ -the one chosen in the article-, which provides values which are only slightly lower than $\hat{\mathrm{V}}_E(\hat{\kappa}_{CU})$ when $n$≥50.

**Acknowledgments**

This research was supported by the Ministry of Science and Innovation (Spain), Grant PID2021-126095NB-I00 funded by MCIN/AEI/10.13039/501100011033 and by “ERDF A way of making Europe”.

**Table 6**

**Results of the 10,000 simulations performed for the two estimated variances of the new Cohen *kappa* coefficient.**

| $K$ | $n$ | $\kappa_C \approx$ | $\hat{V}_E(\hat{\kappa}_{CU})$ | $V(\hat{\kappa}_{CU})$ | $V_A(\hat{\kappa}_{CU})$ | $\bar{\hat{V}}(\hat{\kappa}_{CU})$ | $\bar{\hat{V}}_A(\hat{\kappa}_{CU})$ |
|---|---|---|---|---|---|---|---|
| 2 | 10 | 0.4 | 0.0727 | 0.1235 | 0.1222 | 0.0633 | 0.0627 |
| | | 0.8 | 0.0226 | 0.0187 | 0.0185 | 0.0325 | 0.0322 |
| | 20 | 0.4 | 0.0403 | 0.0382 | 0.0381 | 0.0357 | 0.0356 |
| | | 0.8 | 0.0147 | 0.0139 | 0.0138 | 0.0167 | 0.0167 |
| | 50 | 0.4 | 0.0164 | 0.0162 | 0.0161 | 0.0157 | 0.0157 |
| | | 0.8 | 0.0069 | 0.0064 | 0.0064 | 0.0065 | 0.0065 |
| | 100 | 0.4 | 0.0079 | 0.0079 | 0.0079 | 0.0078 | 0.0078 |
| | | 0.8 | 0.0043 | 0.0041 | 0.0041 | 0.0042 | 0.0042 |
| 3 | 10 | 0.4 | 0.0474 | 0.0385 | 0.0381 | 0.0344 | 0.0341 |
| | | 0.8 | 0.0241 | 0.0126 | 0.0125 | 0.0221 | 0.0219 |
| | 20 | 0.4 | 0.0234 | 0.0218 | 0.0218 | 0.0202 | 0.0201 |
| | | 0.8 | 0.0129 | 0.0107 | 0.0107 | 0.0113 | 0.0113 |
| | 50 | 0.4 | 0.0118 | 0.0115 | 0.0115 | 0.0111 | 0.0111 |
| | | 0.8 | 0.0053 | 0.0049 | 0.0049 | 0.0049 | 0.0049 |
| | 100 | 0.4 | 0.0058 | 0.0058 | 0.0058 | 0.0056 | 0.0056 |
| | | 0.8 | 0.0033 | 0.0032 | 0.0032 | 0.0032 | 0.0032 |
| 5 | 10 | 0.4 | 0.0329 | 0.0260 | 0.0257 | 0.0226 | 0.0224 |
| | | 0.8 | 0.0138 | 0.0088 | 0.0088 | 0.0112 | 0.0111 |
| | 20 | 0.4 | 0.0185 | 0.0172 | 0.0171 | 0.0155 | 0.0155 |
| | | 0.8 | 0.0089 | 0.0074 | 0.0074 | 0.0077 | 0.0077 |
| | 50 | 0.4 | 0.0077 | 0.0075 | 0.0075 | 0.0072 | 0.0072 |
| | | 0.8 | 0.0044 | 0.0038 | 0.0038 | 0.0038 | 0.0038 |
| | 100 | 0.4 | 0.0036 | 0.0035 | 0.0035 | 0.0035 | 0.0035 |
| | | 0.8 | 0.0024 | 0.0023 | 0.0023 | 0.0023 | 0.0023 |

(1) $\hat{V}_E(\hat{\kappa}_{CU})$ is the "exact" variance, or the sample variance of the 10,000 values obtained for $\hat{\kappa}_{CU}$.

(2) $V(\hat{\kappa}_{CU})$ and $V_A(\hat{\kappa}_{CU})$ are the exact asymptotic variances, based on the true values $\{p_{i\bullet}, p_{\bullet j}, p_{ij}, I_e, I_o, \kappa_C\}$.

(3) $\bar{\hat{V}}(\hat{\kappa}_{CU})$ and $\bar{\hat{V}}_A(\hat{\kappa}_{CU})$ are the averages of the 10,000 estimated variances $\hat{V}(\hat{\kappa}_{CU})$ and $\hat{V}_A(\hat{\kappa}_{CU})$, respectively.